\documentclass{amsart}          

\usepackage{amsmath,amsthm}     
\usepackage{amssymb}            
\usepackage{euscript}           
\usepackage{enumerate,calc}
\usepackage[matrix,arrow,curve,frame]{xy}    

\xymatrixcolsep{1.9pc}                          
\xymatrixrowsep{1.9pc}
\newdir{ >}{{}*!/-5pt/\dir{>}}                  

\def\longfib{\DOTSB\relbar\joinrel\twoheadrightarrow}

\newtheorem{thm}[subsection]{Theorem}

\newtheorem{prop}[subsection]{Proposition}

\newtheorem{lemma}[subsection]{Lemma}

\theoremstyle{definition}  

\newtheorem{exercise}[subsection]{Exercise}

\newtheorem{remark}[subsection]{Remark}

  {\end{list}}

\newcommand{\dfn}{\textbf} 

\newcommand{\mdfn}[1]{\dfn{\mathversion{bold}#1}} 

\newcommand{\Smash}             {\wedge}

\newcommand{\dSmash}            {\,{\underline{\wedge}}\,}

\newcommand{\tens}              {\otimes}               
\newcommand{\iso}               {\cong}  

\newcommand{\cat}{\EuScript}    

\newcommand{\cE}{{\cat E}}

\newcommand{\cG}{{\cat G}}
\newcommand{\cH}{{\cat H}}

\newcommand{\Spectra}{{\cat Spectra}}

\newcommand{\Ho}{\text{Ho}}
\newcommand{\ho}{\text{Ho}\,}

\newcommand{\field}[1]  {\mathbb #1} 
\newcommand{\A}         {\field A}
\newcommand{\B}         {\field B}

\newcommand{\Z}         {\field Z}

\newcommand{\W}         {\field W}
\newcommand{\sS}         {\field S}

\DeclareMathOperator*{\Tot}{Tot}

\DeclareMathOperator*{\colim}{colim}

\DeclareMathOperator{\Hom}{Hom}

\DeclareMathOperator{\sk}{sk}

\DeclareMathOperator{\F}{{\mathcal F}}

\DeclareMathOperator{\dF}{\underline{\mathcal F}}

\newcommand{\ra}{\rightarrow}                   
\newcommand{\lra}{\longrightarrow}              
\newcommand{\llra}[1]{\stackrel{#1}{\lra}}      

\newcommand{\cof}{\rightarrowtail}              
\newcommand{\cofib}{\rightarrowtail}              
\newcommand{\trfib}{\stackrel{\sim}{\longfib}}

\newcommand{\inc}{\hookrightarrow}              
\newcommand{\dbra}{\rightrightarrows}           

\newcommand{\blank}{-}                          
\newcommand{\id}{id}                            

\newcommand{\bd}{\partial}

\newcommand{\Loop}{\Omega}

\newcommand{\tW}{\widetilde{W}}
\newcommand{\tC}{\widetilde{C}}

\newcommand{\hW}{\widehat{W}}

\newcommand{\tcE}{\tilde{\cE}}
\newcommand{\tH}{\tilde{H}}
\newcommand{\Si}{\Sigma^\infty}

\newcommand{\shf}{{\text{shf}}}
\newcommand{\sng}{{\text{sing}}}

\newcommand{\HC}{HC}
\newcommand{\trun}[2]{\tau_{{#1}\leq{#2}}}

\newcommand{\vf}{}

\newcommand{\pair}{_{\perp}}

\numberwithin{equation}{section}

\begin{document}

\title{Multiplicative structures on homotopy spectral sequences II}

\author{Daniel Dugger}
\address{Department of Mathematics\\University of Oregon\\ Eugene, OR 97403} 

\email{ddugger@math.uoregon.edu}


\maketitle

\tableofcontents

\section{Introduction}


This short paper is a companion to \cite{multa}.  Here the main
results of that paper are used to establish multiplicative structures
on a few standard spectral sequences.  The applications consist of (a)
applying \cite[Theorem 6.1]{multa} to obtain a pairing of spectral
sequences, and (b) identifying the pairing on the $E_1$- or $E_2$-term
with something familiar, like a pairing of singular cohomology groups.
Most of the arguments are straightforward, but there are subtleties
that appear from time to time.

Originally the aim was just to record a careful treatment of pairings
on Postnikov/Whitehead towers, but in the end other examples have been
included because it made sense to do so.  These examples are just the
ones that I personally have needed to use at some point over the
years, and so of course it is a very limited selection. 

\medskip

In this paper all the notation and conventions of \cite{multa} remain
in force.  In particular, the reader is referred to \cite[Appendix
C]{multa} for our standing assumptions about the category of spaces
and spectra, and for basic results about signs for boundary maps.  The
symbol $\dSmash$ denotes the derived functor of $\Smash$, and
$W_\perp$ denotes an `augmented tower' as in \cite[Section
6]{multa}.  $\Ho(\blank,\blank)$ denotes maps in the homotopy category
of $\Spectra$.
If $A$ is a pointed space we
will write $\F(A,X)$ as an abbreviation for $\F(\Si A,X)$.
Finally, the phrase `globally isomorphic' is often used in the
identification of $E_2$-terms of spectral sequences.  It is explained
in Remark~\ref{re:global}.

\vf

\section{Sign conventions in singular cohomology}
\label{se:prelim}

If $E$ is a ring spectrum and $X$ is a space, then the
Atiyah-Hirzebruch spectral sequence $E_2^{p,q}=H^p(X;E^q) \Rightarrow
E^{p+q}(X)$ is multiplicative.  The naive guess about what
this means is that 
there is an isomorphism of bi-graded rings
\[ 
\oplus_{p,q} E_2^{p,q}(X) \iso \oplus_{p,q} H^p(X;E^q), 
\]
where the products on the right-hand-side are the usual ones
\[ \mu\colon H^p(X;E^q)\tens H^s(X;E^t) \ra H^{p+s}(X;E^{q+t})
\]
induced by the pairings $E^q\tens E^t \ra E^{q+t}$.  Unfortunately,
this statement is just not true in general---one has to add an
appropriate sign into the definition of $\mu$, and these signs cannot
be made to go away.  In order to keep track of such signs in a simple
way, it's useful to re-evaluate the `standard' conventions about
singular cohomology.  I'm grateful to Jim McClure for conversations
about these sign issues. 

\bigskip

Let $X$ be a CW-complex and $C_*(X)$ be its associated cellular chain
complex.  In most algebraic topology textbooks the corresponding
cellular cochain complex is defined by $C^p(X)=\Hom(C_p(X),\Z)$ and
\[ (\delta \alpha)(c)=\alpha(\bd c), \quad \text{for any $\alpha\in C^p(X)$}.
\]
The cup-product of a $p$-cochain $\alpha$ and a $q$-cochain $\beta$ is
defined by the formula
\[ (\alpha\cup\beta)(c\tens d)=\alpha(c) \cdot \beta(d) 
\]
where $c$ is a $p$-chain and $d$ is a $q$-chain.  (Note that we have
written the above formula as if it were an external
cup-product, so we technically need to throw in a diagonal map
somewhere---we omit this to simplify the typography).   
Both of these formulas obviously violate the Koszul sign rule: we will
abandon them and instead define
\begin{equation}
\label{eq:signs}
 (\delta \alpha)(c) = -(-1)^p \alpha(\bd c) \qquad\text{and}\qquad
(\alpha\cup\beta)(c\tens d)=(-1)^{qp}\alpha(c) \cdot \beta(d).
\end{equation}
The first equation may seem to have an unexpected minus sign, but here
is the explanation.  Recall that if $A_*$ and $B_*$ are chain
complexes then there is an associated chain complex $\Hom(A,B)$.  Our
definition of $\delta$ corresponds to the differential on the chain
complex $\Hom(C_*(X),\Z[0])$, where $\Z[0]$ is the complex with  $\Z$
concentrated in dimension $0$.

The sign conventions from (\ref{eq:signs}) appear in \cite{Do}.
We'll of course use these same conventions for cohomology with
coefficients, external cup products, and any similar construction we
encounter.

\begin{exercise}
Check that $\delta$ is a derivation with respect to the cup-product,
and that the dga $C^*(X;\Z)$ defined via our new
formulas is isomorphic to the dga $C^*_{\text{classical}}(X;\Z)$
defined via the old formulas.  In particular, our singular cohomology
ring $H^*(X;\Z)$ is isomorphic to the classical one.
\end{exercise} 

\subsection{Cohomology with graded coefficients}
\label{se:gcup}
If $A_*$ is a graded ring we next want to define the {\it singular
cohomology ring with graded coefficients} $H^*_{grd}(X;A)$, making use
of the natural sign conventions.  It would be nice to just use
the internal hom for chain complexes $\Hom(C_*(X),A)$, where $A$ is
interpreted as having zero differential, but unfortunately this might
give us infinite products in places we don't really want them.  
Instead we'll consider a certain subcomplex.  We
set $C^{p,q}(X;A)=\Hom(C_p(X);A_q)$ and
$C^n_{grd}(X;A)=\oplus_{p-q=n} C^{p,q}(X;A)$---that is, elements of
$C^{p,q}(X;A)$ are regarded as having total degree $p-q$.  For $\alpha
\in C^n_{grd}(X;A)$ we define $\delta\alpha$ by the formula
\begin{equation}
\label{eq:gsign}
 (\delta \alpha)(c)=-(-1)^n \alpha(\bd c).
\end{equation}
The homology of this complex will be denoted $H^*_{grd}(X;A)$; it has
a natural direct sum decomposition into groups $H^{p,q}(X;A)$.

The {\it graded cup-product\/} will be defined on the chain complex
$C_{grd}^*(X;A)$ as follows: if $\alpha \in C^{p,q}(X;A)$ and $\beta \in
C^{s,t}(X;A)$ then
\[ (\alpha\cup \beta)(c\tens
d)=(-1)^{(s-t)p}\alpha(c)\cdot\beta(d) 
\]
(where $c$ is a $p$-chain and $d$ is a $q$-chain).  The sign is again
just the one dictated by the usual Koszul convention, and $\delta$
becomes a derivation with respect to this product.  

If $C_*\tens D_* \ra E_*$ is a pairing of graded abelian groups, one
also has an external graded cup-product $H^*_{grd}(X;C) \tens
H^*_{grd}(Y;D)\ra H^*_{grd}(X\times Y;E)$ defined in a similar
fashion.  It is this product which arises naturally in pairings of
spectral sequences.

\begin{exercise}
\label{ex:signs}
Construct a bi-graded family of isomorphisms
$\eta_{p,q}\colon H^{p,q}(Z;A) \ra H^p(Z;A_q)$, natural in both $Z$
and $A$, which makes the diagrams
\[ \xymatrix{
H^{p,q}(X;C) \tens H^{s,t}(Y;D) \ar[d]_{\eta\tens\eta}\ar[r] & 
      H^{p+s,q+t}(X\times Y;E) \ar[d]_{\eta} \\
H^p(X;C_q) \tens H^s(Y;D_t) \ar[r] & H^{p+s}(X\times Y;E_{q+t})
}
\]
commute up to the sign $(-1)^{sq}$. 
Here $H^n(Z;A_m)$ denotes singular cohomology with coefficients in
$A_m$ as defined via the formulas in (\ref{eq:signs}), and the bottom
map is the cup-product pairing associated to $C_q\tens D_t \ra
E_{q+t}$ (again with the signs from (\ref{eq:signs})).  
\end{exercise}

\begin{exercise}
\label{ex:signs2}
Repeat the above exercise, but this time show that the isomorphisms
$\eta_{p,q}$ can be chosen to make the squares commute up to the sign
$(-1)^{pt}$.  Convince yourself that it is not possible to choose the
$\eta_{p,q}$'s so that the squares commute on the nose.  
\end{exercise}


\vf

\section{Spectral sequences for filtered spaces}
\label{se:cell}

In this section we treat the Atiyah-Hirzebruch spectral sequence, the
Serre spectral sequence, and spectral sequences coming from geometric
realizations.  Some other references for the former are \cite{K},
\cite[Appendix B]{GM}, \cite{V}.  For the Serre spectral sequence see
\cite{K}, \cite[Chap. 5]{Mc}, \cite[Chap. 9.4]{Sp}, and \cite[XIII.8]{Wh}.

\subsection{Generalities}
\label{se:filt}
Suppose given a sequence of cofibrations $\emptyset \cofib A_0 \cofib
A_1 \cofib \cdots$ and let $A$ denote the colimit.  If $\emptyset
\cofib B_0 \cofib B_1 \cofib \cdots \cofib B$ is another sequence of
cofibrations, we may form the product sequence whose $n$th term is
\[ (A\times B)_n = \bigcup_{i+j=n} (A_i\times B_j).
\]
This is a sequence of cofibrations whose colimit is $A\times B$.
Given a fibrant spectrum $\cE$, one can look at the induced tower
\[  \cdots \ra \F(A_{2+},\cE) \ra \F(A_{1+},\cE) \ra \F(A_{0+},\cE)
\]
and identify the homotopy fibers as $\F(A_n/A_{n-1},\cE)$.  This is a
lim-tower rather than a colim-tower, and is not convenient for seeing
multiplicative structures; one doesn't have reasonable pairings
$\F(A_{k},\cE)\Smash \F(B_{n},\cE) \ra \F((A\times B)_{{k+n}},\cE)$,
for instance.  Instead we have to use a slightly different tower.

The cofiber sequences $A_{n}/A_{n-1} \inc A/A_{n-1} \ra A/A_{n} $
induce rigid homotopy fiber sequences $\F(A/A_{n},\cE) \ra
\F(A/A_{n-1},\cE)\ra \F(A_{n}/A_{n-1},\cE)$.  We define an augmented
colim-tower by setting $\W(A,\cE)_n=\F(A/A_{n-1},\cE)$ and
$\B(A,\cE)_n=\F(A_{n}/A_{n-1},\cE)$.  The associated spectral sequence
$E_*(A,\cE)$ might be called the \mdfn{$\cE$-spectral sequence for the
filtered space $A$}.

\begin{exercise}
Verify that the tower $\W(A,\cE)$ is weakly equivalent to the tower
$\Omega \F(A_*,\cE)$ via a canonical zig-zag of towers.  So the
homotopy spectral sequences can be identified. 
\end{exercise}

Now assume that $\cE$ had a multiplication $\cE \Smash \cE \ra \cE$.
Then for any two pointed spaces $X$ and $Y$ we have the map 
\[\F(X,\cE)
\Smash \F(Y,\cE) \ra \F(X\Smash Y,\cE\Smash \cE) \ra \F(X\Smash Y,\cE).
\]
Using this, the obvious maps of spaces
\begin{align*}
(A\times B)/(A\times
B)_{q+t-1} &\ra
A/A_{q-1}\Smash B/B_{t-1},
\quad\text{and} \\
(A\times B)_{q+t}/(A\times B)_{q+t-1}&\ra
A_q/A_{q-1}\Smash B_t/B_{t-1}
\end{align*}
give
pairings $\W(A,\cE) \Smash \W(B,\cE) \ra \W(A\times B,\cE)$ and
$\B(A,\cE) \Smash \B(B,\cE) \ra \B(A\times B,\cE)$ which are
compatible with the maps in the towers.  So
\cite[Thm 6.1]{multa} gives us a pairing of spectral
sequences $E_*(A,\cE) \tens E_*(B,\cE) \ra E_*(A\times B,\cE)$.  
This is the `formal' part of the construction.

\subsection{The Atiyah-Hirzebruch spectral sequence}
\label{se:ah}
Now we specialize to where $A$ and $B$ are CW-complexes which are
filtered by their skeleta.  In this case we can identify the $E_1$-
and $E_2$-terms, and we will need to be very explicit about how we do
this.  For convenience we take $A$ and $B$ to be {\it labelled
CW-complexes\/}, meaning that they come with a chosen indexing of
their cells.  Let $I_q$ be the indexing set for the $q$-cells in $A$,
and let $C_q(A)=\oplus_{\sigma\in I_q} \Z$.

Recall that
\[ E_1^{p,q}=\pi_p\F(A_q/A_{q-1},\cE)\cong \Ho(S^p\Smash
A_q/A_{q-1},\cE).
\]  
An element
$\sigma\in I_q$ specifies a map $S^q \ra A_q/A_{q-1}$ which we will
also call $\sigma$.  
Given an element $f\in E_1^{p,q}=\Ho(S^p\Smash A_q/A_{q-1},\cE)$,
restricting to each $\sigma$ specifies an element in $\Ho(S^p\Smash
S^q,\cE)$.  We therefore get a cochain in
$\Hom(C_q(A),\pi_{p+q}\cE)=C^{q,p+q}_{grd}(A;\cE_*)$, and we'll choose
this assignment for our isomorphism
\[ E_1^{p,q} \iso C^{q,p+q}_{grd}(A;\cE_*).
\]
Note that this isomorphism is completely natural with respect to maps
of labelled CW-complexes.

We claim that the $d_1$-differential corresponds under this
isomorphism to the differential on $C^*_{grd}$ defined in
Section~\ref{se:prelim}.  By naturality (applied a couple of times) it
suffices to check this when $A$ is the CW-complex with $A_{q-1}=*$,
$A_{q}=S^{q}$, and $A_k=D^{q+1}$ for $k\geq q+1$.  In this case
our $d_1$ is the boundary map in the
long exact homotopy sequence of $\F(S^{q+1},\cE) \ra \F(D^{q+1},\cE)
\ra \F(S^q,\cE)$, which takes the form
\[ \ho\!(S^p\Smash S^q,\cE) = \ho\!(S^p,\F(S^q,\cE)) \llra{\bd}
\ho(S^{p-1},\F(S^{q+1},\cE)) = \ho\!(S^{p-1}\Smash S^{q+1},\cE).
\]
We know from \cite[C.6(d)]{multa} that the composite is $(-1)^{p-1}$
times the canonical map.  Via our identification with cochains, we are
looking at a map $C^q(A;\pi_{p+q}\cE) \ra C^{q+1}(A;\pi_{p+q}\cE)$,
and the sign $(-1)^{p-1}$ is precisely the one for the coboundary
$\delta$ defined in (\ref{eq:gsign}).

In a moment we will identify the pairing on $E_2$-terms, but before
that we make a brief remark on the case $A=B$.  The diagonal map $A\ra
A\times A$ is homotopic to a map $\Delta'$ which preserves the
cellular filtration, and so $\Delta'$ induces a map of towers
$\W_\perp(A\times A,\cE) \ra \W_\perp(A,\cE)$.  Composing this with
our above pairing gives $\W_\perp(A,\cE)\Smash \W_\perp(A,\cE) \ra
\W_\perp(A,\cE)$, and so we get a multiplicative structure on the
spectral sequence $E_*(A,\cE)$.

\begin{thm}[Multiplicativity of the Atiyah-Hirzebruch spectral sequence]
\label{th:cellular}
There is a natural pairing of spectral sequences $E_*(A,\cE) \tens
E_*(B,\cE) \ra E_*(A\times B,\cE)$ together with natural isomorphisms
$\oplus_{p,q} E_2^{p,q}(?,\cE) \iso \oplus_{p,q}H^q(?,\cE^{-p-q})$
(for $?=A,B,A\times B$) which make the diagrams
\[ \xymatrix{
E_2^{p,q}(A,\cE) \tens E_2^{s,t}(B,\cE) \ar[r]\ar[d] 
& E_2^{p+s,q+t}(A\times B,\cE) \ar[d] \\
H^q(A;\cE^{-p-q}) \tens H^t(B;\cE^{-s-t}) \ar[r] 
& H^{q+t}(A\times B;\cE^{-p-q-s-t})
}
\]
commute, where the bottom map is the graded cup product from
Section~\ref{se:prelim}.

In the diagonal case, there is a natural isomorphism of rings
$\oplus_{p,q}E_2^{p,q}(A,\cE) \iso \oplus_{p,q}H^q(A;\cE^{-p-q})$,
where the latter is again given the graded cup product.
\end{thm}

\begin{remark}
\label{re:global}
Rather than repeat the above statement for every multiplicative
spectral sequence we come across, we'll just say that the $E_2$-term
is \dfn{globally isomorphic} to the graded cup product (that is, they
are naturally isomorphic as pairings of bigraded abelian groups).
\end{remark}

\begin{proof}
We have done everything except identify the product.  
The pairing on $E_1$-terms is the map
\[ \pi_p \F(A_q/A_{q-1},\cE) \Smash \pi_s\F(B_t/B_{t-1},\cE) \ra
\pi_{p+s} \F(A_q/A_{q-1}\Smash B_t/B_{t-1},\cE).
\]
Recall that this sends $\alpha\colon S^p\Smash A_q/A_{q-1} \ra \cE$
and $\beta\colon S^s\Smash B_t/B_{t-1} \ra \cE$ to the composite
\[ \alpha\beta\colon 
S^p\Smash S^s \Smash A_q/A_{q-1}\Smash B_t/B_{t-1} \ra
S^p \Smash A_q/A_{q-1} \Smash S^s \Smash B_t/B_{t-1} \ra \cE \Smash
\cE \ra \cE.
\]
Choosing a $q$-cell $\sigma$ of $A$ yields a map $S^q \ra
A_q/A_{q-1}$,
and a $t$-cell $\theta$ of $B$ gives a map $S^t \ra B_t/B_{t-1}$.  
Under our identification with cochains, the `value' of $\alpha\beta$
on the cell $\sigma\Smash \theta$ is the restriction of $\alpha\beta$
to $S^p\Smash S^s \Smash S^q \Smash S^t$.  

If, on the other hand, we compute $\alpha(\sigma)\cdot \beta(\theta)$
in the ring $\pi_*\cE$, we get the composite
\[ [S^p\Smash S^q] \Smash [S^s \Smash S^t] \ra [S^p\Smash A_{q}/A_{q-1}]
\Smash [S^s \Smash B_t/B_{t-1}] \ra \cE \Smash \cE \ra \cE.
\]
By inspection, this differs from $(\alpha\beta)(\sigma\Smash \theta)$
by the sign $(-1)^{sq}$, which is the same sign that was used in
defining the graded cup product from section~\ref{se:gcup} (remember
that under our isomorphism $\alpha$ lies in $C^{q,p+q}(A;\cE_*)$ and
$\beta$ lies in $C^{t,s+t}(B;\cE_*)$).
\end{proof}

\begin{remark}
In the square from the statement of the theorem, the bottom map is
$(-1)^{t(p+q)}$ times the `un-graded' cup product on cohomology
induced by the pairing $\cE^{-p-q} \tens \cE^{-s-t} \ra
\cE^{-p-q-s-t}$.  This follows from Exercise~\ref{ex:signs}.  The
signs are easy to remember, because they follow the usual conventions:
The index `$t$' is commuted across the index `$-(p+q)$', and as a
result the sign $(-1)^{-t(p+q)}$ is picked up (the minus sign can of
course be left off the exponent).  Note that for most of the familiar
cohomology theories, like $K$-theory or complex cobordism, the signs
end up being irrelevant because the coefficient groups are
concentrated in even dimensions.
\end{remark}

\subsection{The Serre spectral sequence}
%
%
Let $p\colon X \ra B$ be a fibration with fiber $F$, where $B$ is a
pointed, connected CW-complex.  Let $B_0 \subseteq B_1 \subseteq
\cdots$ be the skeletal filtration of $B$, and define $X_i=p^{-1}B_i$.
We'll assume that the inclusions $X_i \inc X_{i+1}$ are cofibrations
between cofibrant objects, and consider the augmented tower
$\W_n=\F(X/X_{n-1},H\Z)$, $\B_n=\F(X_n/X_{n-1},H\Z)$.  The associated
homotopy spectral sequence $E_*(X)$ is the \dfn{Serre spectral
sequence} for the fibration.

It is easy to see that there is a natural identification 
\[ X_n/X_{n-1}
\iso \bigvee_\alpha \Bigl [p^{-1}e^n_\alpha /p^{-1} \bd(e^n_\alpha)
\Bigr ],
\]
where the wedge ranges over the $n$-cells $e^n_\alpha$ of $B$.  The
interior of a cell $e^n$ is just the interior of $D^n$, so we can take
a closed disk around the origin with radius $\frac{1}{2}$---call this
smaller disk $U$.  Then $\tH^*(p^{-1}e^n/p^{-1}\bd(e^n))$ may be
canonically identified with $\tH^*(p^{-1}U/p^{-1}\bd U)$, and we are
better off than before because $\bd U$ is actually a sphere (rather
than just the image of one).  The
diagram
\[ \xymatrix{ p^{-1}(0) \ar[r]^\sim
\ar@{ >->}[d]_\sim &p^{-1}U \ar@{->>}[d] \\
             U\times p^{-1}(0) \ar@{.>}[ur]^\lambda\ar@{->>}[r] & U
}
\]
has a lifting as shown, and this lifting will be a weak equivalence.
It restricts to a weak equivalence $\bd U\times p^{-1}(0) \ra
p^{-1}(\bd U)$ (because this is a map of fibrations over $\bd U$, and
it is a weak equivalence on all fibers).  Therefore we have the
diagram
\[ \xymatrix{
{*}  & p^{-1}(\bd U) \ar[l]\ar@{ >->}[r] & p^{-1}(U) \\
{*} \ar@{=}[u]  & \bd U \times p^{-1}(0) \ar[u]_\sim\ar[l]\ar@{ >->}[r] 
& U \times p^{-1}(0) \ar[u]_\sim
}
\]
and this necessarily induces a weak equivalence on the pushouts.  
In this way we get an identification 
\[
\tH^k(p^{-1}U/p^{-1}(\bd U)) \iso\!
\tH^k([U/\bd U] \Smash p^{-1}(0)_+)
\iso \tH^k(S^n \Smash p^{-1}(0)_+)\iso\! H^{k-n}(p^{-1}(0)).
\]
Of course the first isomorphism depended on the lifting $\lambda$, and
so is not canonical.

We refer to \cite[Chap. 5]{Mc} for a detailed discussion of local
coefficient systems and their use in this particular context.  But
once the right definitions are in place the argument we gave for the
Atiyah-Hirzebruch spectral sequence in the last section adapts
verbatim to naturally identify the $(E_1,d_1)$-complex as
\[ E_1^{p,q}\iso \pi_p \F(X_q/X_{q-1},H\Z) \iso C^q(B; \cH^{p+q}(F))
\]
where $\cH^*(F)$ denotes the appropriate system of coefficients.  The
differential on the cochain complex is still the one from
section~\ref{se:prelim}, appropriate for cellular cohomology with
graded coefficients.  

Now suppose $X'\ra B'$ is another fibration satisfying the same basic
assumptions as $X\ra B$.  We give $B\times B'$ the product cellular
filtration, and then pull it back to get a corresponding filtration of
$X\times X'$.  This coincides with the product of the filtrations on
$X$ and $X'$, and so we get a pairing of spectral sequences by the
discussion in section~\ref{se:filt}.  The identification of the
pairing with the graded cup product again follows exactly as for the
Atiyah-Hirzebruch spectral sequence.

\begin{thm}[Multiplicativity of the Serre spectral sequence]
There is a natural pairing of Serre spectral sequences $E_*(X) \tens
E_*(X') \ra E_*(X\times X')$ such that the pairing of $E_2$-terms is
globally isomorphic to the graded cup product on singular cohomology
with local coefficients.
\end{thm}

\subsection{Spectral sequences for simplicial objects}
\label{se:simpspace}

Filtered spaces also arise in the context of geometric realizations.
Let $X_*$ be a Reedy cofibrant simplicial space, in which case the
skeletal filtration of the realization $|X|$ is a sequence
of cofibrations.  There is a resulting tower of rigid homotopy cofiber
sequences with $W_q(X_*,\cE)=\F(|X|/\sk_{q-1}|X|,\cE)$ and
$B_q(X_*,\cE)=\F(\sk_q|X|/\sk_{q-1}|X|,\cE)$.

If $Y_*$ is another Reedy cofibrant simplicial space, we can equip
$|X|\times |Y|$ with the product filtration.  We also have the product
simplicial space $X_*\times Y_*$, equipped with its skeletal
filtration.  There is a natural map $\eta\colon |X\times Y| \ra
|X|\times |Y|$, and this is actually an isomorphism (using adjointness
arguments one reduces to the case where $X$ and $Y$ are the simplicial
sets $\Delta^m$ and $\Delta^n$, where it is (T2) from \cite[Appendix
C]{multa}).  Unfortunately $\eta$ does {\it not\/} preserve the
filtrations, as can be seen by taking $X$ and $Y$ both to be the
simplicial set $\Delta^1$ (regarded as a discrete simplicial space).
The product filtration on $|\Delta^1|\times |\Delta^1|$ is smaller
than the skeletal filtration coming from $|\Delta^1\times\Delta^1|$.

The formal machinery of section~\ref{se:filt} gives a pairing from
$E_*(|X|,\cE)$ and $E_*(|Y|,\cE)$ to the spectral sequence for the
product filtration on $|X|\times |Y|$---let's call this
$E_*(|X|\times|Y|,\cE)$.  Often one would like to have a pairing into
$E_*(|X\times Y|,\cE)$, but this doesn't seem to follow from our basic
results.  Here are two ways around this.  One can replace
$E_*(|X|,\cE)$ with the homotopy spectral sequence for the
cosimplicial spectrum $[n]\mapsto \F(X_n,\cE)$, and similarly for
$E_*(|Y|,\cE)$ and $E_*(|X\times Y|,\cE)$.  The paper \cite{BK}
proves that if $M^* \Smash N^* \ra Q^*$ is a (level-wise) pairing of
cosimplicial spaces, then there is an associated pairing of spectral
sequences---this gives us what we wanted.  Having not checked the
details in \cite{BK}, I can say nothing more about this approach;
their results clearly depend on more than the formal theorems of
\cite{multa}, but I couldn't tell from their paper exactly what
the important ingredient is.

Here is another approach which sometimes works.  While $\eta
\colon |X\times Y| \ra |X|\times|Y|$ does not preserve the
filtrations, $\eta^{-1}$ {\it is\/} filtration-preserving (by
functoriality and adjointness arguments it suffices to check this when
$X$ and $Y$ are the simplicial sets $\Delta^n$ and $\Delta^m$).  So
$\eta^{-1}$ induces a map of spectral sequences $E_*(|X\times Y|,\cE)\ra
E_*(|X|\times |Y|,\cE)$.  We have the following:

\begin{prop}
\label{pr:simp}
If $X_*$ and $Y_*$ are simplicial sets, the natural map of spectral
sequences $E_*(|X\times Y|,\cE)\ra E_*(|X|\times |Y|,\cE)$ is an
isomorphism on $E_2$-terms.
\end{prop}

\begin{proof}
This follows from the work in section \ref{se:ah}, since it
identifies both $E_1$-terms as cellular chain complexes computing
$H^*(|X\times Y|,\cE^*)$, but for different CW-decompositions.
\end{proof}

It follows that when $X_*$ and $Y_*$ are simplicial sets we get our
desired pairing $E_*(|X|,\cE)\tens E_*(|Y|,\cE) \ra E_*(|X\times
Y|,\cE)$ from the $E_2$-terms onward.  This observation will be used
in section~\ref{se:open}.

\begin{exercise}
Is Proposition~\ref{pr:simp} true for simplicial spaces?  I haven't
worked out the answer to this.
\end{exercise}

\vf 


\section{The Postnikov/Whitehead spectral sequence}
\label{se:postnikov}

For each spectrum $E$ and each $n\in \Z$, let $P_nE$ denote the $n$th
Postnikov section of $E$; this is a spectrum obtained from $E$ by
attaching cells to kill off all homotopy groups from dimension $n+1$
and up.  The construction can be set up so that if $E$ is fibrant then all
the $P_nE$ are also fibrant, and there are natural maps
$E\ra P_n E$ and $P_n E \ra P_{n-1}E$ making the obvious triangle
commute.  So we have a tower of fibrant spectra
\[ \cdots \ra P_2 E \ra P_1 E \ra P_0 E \ra \cdots 
\]
and the homotopy cofiber of $P_{n+1} E \ra P_{n}E$ is an Eilenberg-MacLane
spectrum of type $\Sigma^{n+2} H(\pi_{n+1}E)$.

If $A$ is a cofibrant, pointed space, we can map $\Si A$ into this
tower and thereby get a tower of function spectra
\[ \cdots \ra \F(A,P_2E) \ra \F(A,P_1 E) \ra \F(A,P_0 E) \ra \cdots 
\]
The homotopy cofiber of $\F(A,P_{n+1}E) \ra \F(A,P_{n}E)$ is weakly
equivalent to $\dF(A,\Sigma^{n+2} H(\pi_{n+1}E))$, and the resulting
homotopy spectral sequence has
\[ E_1^{p,q}\cong \pi_p \dF(A,\Sigma^{q+2} H(\pi_{q+1} E))\cong 
H^{q-p+2}(A;\pi_{q+1} E).
\]
The spectral sequence abuts to $\pi_{p-1} \dF(A,E) =
\tilde{E}^{1-p}(A)$.  This turns out to be another construction of the
Atiyah-Hirzebruch spectral sequence---see \cite[Appendix]{GM} for some
information about how the two spectral sequences are related.

Assume that $E$, $F$, and $G$ are fibrant spectra, and that there is a
pairing $E\Smash F \ra G$.  There do not exist reasonable pairings
$P_nE\Smash P_m F \ra P_{n+m}G$, and so Postnikov towers are not
convenient for seeing multiplicative structures on spectral sequences.
This is related to the Postnikov tower being a lim-tower rather than a
colim-tower.  Instead we will use the `reverse' of the Postnikov
tower, sometimes called the Whitehead tower.  If $W_nE$ denotes the
homotopy fiber of $E \ra P_{n-1}E$, then there are natural maps $W_n E
\ra W_{n-1} E$ and so we get a new tower.  The homotopy cofiber of
$W_{n+1} E \ra W_n E$ is weakly equivalent to $\Sigma^n H(\pi_n E)$.
We will modify these towers in an attempt to produce a
pairing $W_* E \Smash W_* F \ra W_* G$.

To explain the idea, let's forget about cofibrancy/fibrancy issues
for just a moment.  Consider the following maps:
\[ \xymatrix{
  W_m E \Smash W_n F \ar@{.>}[d]_{\lambda}\ar[r]
&E\Smash F\ar[d] \\
W_{m+n}G \ar[r]^-p & G \ar[r]^-j & P_{m+n-1}G.
}
\]
The horizontal row is a homotopy fiber sequence.  The spectrum $W_m E
\Smash W_n F$ is $(m+n-1)$-connected, and so the composite $W_mE
\Smash W_n F \ra P_{m+n-1}G$ is null-homotopic.  Choosing a
null-homotopy lets us construct a lifting $\lambda\colon W_m E\Smash
W_n F \ra W_{m+n}G$.  If we had two different liftings $\lambda$ and
$\lambda'$, their difference would lift to a map $W_mE \Smash W_n F
\ra \Loop P_{m+n-1}G$ and so would be null-homotopic (again, because
the domain is $(m+n-1)$-connected).  So the lift
$\lambda$ is unique up to homotopy.

The situation, then, is that we can produce pairings $W_m E\Smash W_n
F \ra W_{m+n}G$, but so far they don't necessarily commute with the
structure maps in the towers.  They certainly commute up to
homotopy---this follows from the `uniqueness' considerations in the
above paragraph--- but we need them to commute on the nose.  By using
obstruction theory we will be able to alter these maps so that the
relevant diagrams do indeed commute.  The argument proceeds in a few
steps.

\begin{lemma}
\label{le:post2}
For each fibrant spectrum $E$ there is a natural tower of rigid
homotopy cofiber sequences $(\tW_*E,\tC_*E)$ such that every $\tW_n E$
and $\tC_n E$ is cofibrant-fibrant, together with a natural zig-zag of
weak equivalences from $\tW_*E$ to $W_*E$.
\end{lemma}

\begin{proof}
First take $W_*E$ and apply a cofibrant-replacement functor $Q$ to
all the levels: this produces $QW_*E$, a tower of cofibrant spectra.
Then perform the telescope construction from \cite[B.4]{multa} to
get a tower of cofibrations between cofibrant objects $TW_*E$ and a
weak equivalence $TW_*E \ra QW_*E$.  Let $C_nE$ denote the cofiber of
$TW_{n+1}E\ra TW_nE$.  Finally, let $F$ be the fibrant-replacement
functor for $\Spectra$ such that $F(*)=*$ given in \cite[C.3(c)]{multa}.
Applying $F$ to the rigid tower $(TWE,CE)$ gives a new rigid tower
which has the desired properties.
\end{proof}

At this point we have towers where everything is cofibrant-fibrant, so
the argument we have already explained will construct maps $\tW_mE
\Smash \tW_n F \ra \tW_{m+n}G$ which commute up to homotopy with the
maps in the towers.  By considering the diagram
\[ \xymatrix{
\tW_mE \Smash \tW_n F \ar[r]\ar[d] & \tW_{m+n}G \ar[d] \\
\tC_m E\Smash \tC_n F   & \tC_{m+n}G
}
\]
one can see that there is a unique homotopy class $\tC_m E \Smash
\tC_n F \ra \tC_{m+n}G$ which makes the square commute.  This is
because $\tC_{m+n}G$ is an Eilenberg-MacLane spectrum of type
$\Sigma^{n+m}H(\pi_{n+m}G)$, and $\tW_mE \Smash \tW_n F \ra \tC_mE
\Smash \tC_nF$ induces an isomorphism on the corresponding cohomology
group (since both the domain and codomain are $(m+n-1)$-connected).
So at this point we have produced a {\it homotopy-pairing\/} $(\tW
E,\tC E) \dSmash (\tW F,\tC F) \ra (\tW G, \tC G)$ (see
\cite[6.3]{multa}).  We will prove:

\begin{prop}
\label{pr:post1}
The homotopy-pairing $(\tW E,\tC E) \dSmash (\tW F,\tC F) \ra (\tW G,
 \tC G)$ is locally realizable.
\end{prop}

The following lemma encapsulates the basic facts we will need.
The proof will be left to the reader.

\begin{lemma}[Obstruction theory]
Suppose that $X\ra Y$ is a fibration of spectra which induces
isomorphisms on $\pi_k$ for $k\geq n$.  Let $A\cof B$ be a cofibration
which induces isomorphism on $\pi_{k}$ for $k<n$.  Then any diagram
\[ \xymatrix{ A \ar[r]\ar@{ >->}[d] & X \ar@{->>}[d] \\
     B \ar[r]\ar@{.>}[ur] & Y}
\]
has a lifting as shown.
\end{lemma}

\begin{proof}[Proof of Proposition~\ref{pr:post1}]
First we truncate the towers, and we might as well assume we are
dealing with truncations $\trun{0}{k}(\tW E,\tC E)$ and
$\trun{0}{l}(\tW F,\tC F)$ because the argument will be the same no
matter what the lower bounds are.  For the rest of the argument we
will only be dealing with these finite towers, and will omit the
$\tau$'s from the notation.

We replace $(\tW E,\tC E)$ and $(\tW F,\tC F)$ by the
equivalent towers $(TW E, CE)$ and $(TW F, CF)$ constructed in the
proof of Lemma~\ref{le:post2}, because these consist of cofibrations
between cofibrant spectra.  It is easy to see that one can also find a
tower $\hW_*G$ consisting of fibrant spectra and {\it fibrations\/},
together with a weak equivalence $\tW_*G \ra \hW_* G$ (remember that
all our towers are finite).  We will construct a pairing of towers
$TW_*E\Smash TW_* F \ra \hW_*G$ which realizes the
homotopy-pairing. For $TW_{0}E \Smash TW_{0}F \ra \hW_{0}G$ we choose
any map in the correct homotopy class.  Next consider the diagram
\[ \xymatrix{ &&   \hW_1 G \ar[d] \\
 TW_0 E\Smash TW_1 F \ar[r] & TW_0 E\Smash TW_0 F \ar[r]
   & \hW_0 G.
}
\]
The vertical map is a fibration inducing isomorphisms on $\pi_1$ and
higher, and the spectrum $TW_0E\Smash TW_1 F$ is 0-connected; so there
is a lifting $\mu_{(0,1)}$.  Next, look at the diagram
\[ \xymatrix{
TW_1 E \Smash TW_1 F \ar[r] \ar[d] & 
TW_0 E\Smash TW_1 F \ar[r] &  \hW_1 G \ar[d] \\
TW_1 E\Smash TW_0 F \ar[r] & TW_0 E\Smash TW_0 F \ar[r] &
\hW_0 G.
}
\]
This diagram commutes (because the missing vertical arrow in the
middle may be filled in). The right vertical map is a fibration which
induces isomorphisms on $\pi_1$ and higher.  The left vertical map
induces isomorphisms on $\pi_0$ and lower (because both domain and
range are $0$-connected), and is a cofibration.  So there is a lifting
$\mu_{(1,0)}\colon TW_1 E\Smash TW_0 F \ra \hW_1 G$.

This process may be continued to inductively define $\mu_{(0,2)}$,
$\mu_{(1,1)}$, and $\mu_{(2,0)}$, and then onward from level three.
In this way, we construct the required pairing of towers.  This
pairing agrees with the original homotopy-pairing  because of
the `uniqueness' of the liftings $\lambda$ in our original discussion;
the details are left for the reader. 

At this point we have a pairing of towers, but we need a pairing of
augmented rigid towers.  We can't just take cofibers in $\hW_*G$
because the maps are fibrations, not cofibrations.  So let $Q(\hW_*G)
\trfib \hW_*G$ be the cofibrant approximation guaranteed in
\cite[Lemma B.2]{multa}.  We have a diagram
\[ \xymatrix{ & Q(\hW_*G) \ar@{->>}[d]^\sim \\
    TW E \Smash TW F \ar[r] & \hW_*G
}
\]
and by \cite[Lemma B.3]{multa} the lower left corner is a cofibrant
tower; so there is a lifting.  The tower $Q(\hW G)$ consists of
cofibrations, and so augmenting by the cofibers gives a rigid tower.
The new pairing automatically passes to cofibers to give $(TW E,C
E)\Smash (TW F,C F) \ra (Q(\hW G),CQ(\hW G))$.

Finally, consider
\[ \xymatrix{
& Q(\hW_*G) \ar@{->>}[d]^\sim \\
\tW_*G \ar[r]^\sim & \hW_*G.
}
\]
The tower $\tW_*G$ was cofibrant, so there is a lifting.  This will be
a levelwise equivalence, and therefore induces an equivalence on the
cofibers.  So we get an equivalence of augmented towers
$(\tW_*G,\tC_*G) \ra (Q(\hW G),CQ(\hW G))$.  Thus, we have constructed
the required realization of our homotopy-pairing.
\end{proof}

For each cofibrant space $X$ consider the tower of function spectra
$\F(X_+,\tW E_{\perp})$ (recall that the $\tW_*E$ are all fibrant).
This is a tower of rigid homotopy cofiber sequences, and we will call
the associated spectral sequence $E_*(X,E)$ the \mdfn{Whitehead
spectral sequence for $X$ based on $E$}.  The homotopy-pairing $\tW
E\pair\dSmash \tW F\pair \ra \tW G\pair$ induces a homotopy-pairing on
towers of function spectra, and by \cite[Prop. 6.10]{multa} this is
locally-realizable and so induces a pairing of spectral sequences: for
any cofibrant spaces $X$ and $Y$ we have $E_*(X,E) \tens E_*(Y,F) \ra
E_*(X\times Y,G)$.  What is left is to identify the pairing on
$E_1$-terms with the pairing on singular cohomology (up to the correct
sign).

If $X$ is a spectrum with a single non-vanishing homotopy group in
dimension $m$, there is a unique isomorphism in the homotopy category
$\sS^m \dSmash H(\pi_m X) \ra X$ with the property that
 the composite
\[ \pi_m X \ra \pi_0 H(\pi_m X) \llra{\sigma_l} \pi_m(\sS^m\Smash
H(\pi_m X)) \ra \pi_m X
\]
is the identity map (the first map in the composite is the one
provided by \cite[Section C.7]{multa}).
If $X \Smash Y \ra Z$ is a pairing of spectra
where $X$, $Y$, and $Z$ each have a single non-vanishing homotopy group
in dimensions $m$, $n$ and $m+n$, then the diagram in $\Ho(\Spectra)$
\[ \xymatrix{
\sS^m \dSmash H(\pi_m X) \dSmash \sS^n\Smash H(\pi_n Y) \ar[r]\ar[d] &
\sS^{m+n} \dSmash H(\pi_{m+n} Z) \ar[d] \\
X\Smash Y \ar[r] & Z
}
\]
is commutative.  Here the top map interchanges the $\sS^n$ and the
$H(\pi_m X)$ and then uses the map $H(\pi_m X)\tens H(\pi_n Y) \ra
H(\pi_{m+n}Z)$ induced by the pairing $\pi_m X\tens \pi_n Y \ra
\pi_{m+n}Z$ (cf. \cite[Section C.7]{multa}).  The above observations are
simple calculations in the homotopy category of spectra.

In our situation we have specific isomorphisms $\pi_m \tC_mE \iso
\pi_m E$, and the same for $F$ and $G$.  This is because $\tW_m E \ra
\tC_m E$ induces an isomorphism on $\pi_m$, $\tW_m E$ is connected by
a chosen zig-zag of weak equivalences to $W_m E$, and the map $W_m E
\ra E$ induces an isomorphism on $\pi_m$ as well.  The pairing $\tC_m
E \Smash \tC_n F \ra \tC_{m+n} G$ induces a pairing on homotopy groups
which corresponds to the expected pairing $\pi_m E\tens \pi_n F \ra
\pi_{m+n}G$ under these isomorphisms.  Putting all the above
statements together, we have proven:

\vfill\eject

\begin{lemma}
\label{le:EMunique}
In $\Ho(\Spectra)$ there exist isomorphisms $\sS^n\Smash H(\pi_n E) \ra
\tC_nE$, $\sS^n \Smash H(\pi_n F) \ra \tC_n F$, and $\sS^n\Smash H(\pi_n G)
\ra \tC_n G$ for all $n\in \Z$, such that the diagrams
\[ \xymatrix{
\sS^m\Smash H(\pi_mE) \Smash \sS^n\Smash H(\pi_n F) \ar[r]\ar[d] &
\sS^{m+n} \Smash H(\pi_{m+n}G) \ar[d] \\
\tC_mE \Smash \tC_n F \ar[r] &\tC_{m+n}G
}
\]
all commute (in the homotopy category).
\end{lemma}

The lemma tells us that if $A$ and $B$ are spectra the pairing
\[
\biggl [ \oplus_{p,q}\Ho(\sS^p\dSmash A,\tC_q E) \biggr ] \tens
\biggl [
\oplus_{s,t}\Ho(\sS^s\dSmash B,\tC_t F) 
\biggr ]
\ra
\oplus_{u,v} \Ho(\sS^{u}\dSmash A\dSmash B,\tC_{v} G)
\]
is globally isomorphic to the pairing obtained from the maps
\[\xymatrix{
\Ho(\sS^p\dSmash A,\sS^q\dSmash H(\pi_qE) ) \tens
\Ho(\sS^s\dSmash B,\sS^t\dSmash H(\pi_t F)) \ar[d] \\
 \Ho(\sS^{p+s}\dSmash A\dSmash B,\sS^{q+t}\dSmash H(\pi_{q+t} G)).
}
\]
Now, the left-suspension map gives an isomorphism
\[
\Ho(\sS^p\dSmash
B,\sS^q \dSmash H(\pi_q E)) \iso \Ho(\sS^{p-q}\dSmash B,H(\pi_q E)) =
H^{q-p}(B;\pi_q E),
\]
and similarly for the $F$ and $G$ terms.  This allows us to rewrite
the above pairing as a pairing of singular cohomology groups, but the
suspension maps introduce signs.  For any spectra $M$ and $N$, the 
diagram
\[ \xymatrixcolsep{1.91pc}\xymatrix{
\ho(\sS^{a}\!\dSmash\! A,M) \tens
\ho(\sS^{b}\!\dSmash\! B,N) \ar[d]\ar[r]^-{\sigma_l^q\Smash \sigma_l^t}
&
\ho(\sS^{q+a}\!\dSmash\! A,\sS^q\!\dSmash\! M) \tens
\ho(\sS^{t+b}\!\dSmash\! B,\sS^t\!\dSmash\! N) \ar[d] \\
\ho(\sS^{a+b}\dSmash A\Smash B, M\dSmash N) \ar[r]^-{\sigma_l^{q+t}}
& 
\ho(\sS^{q+t+a+b}\Smash A\Smash B,\sS^{q+t}\Smash M\Smash N)
}
\]
commutes up to the sign $(-1)^{ta}$ (one compares the string `$qatb$'
to the string `$qtab$' and sees that the $t$ and $a$ must be commuted).
Taking $A=\Si X_+$ and $B=\Si Y_+$, we now conclude:

\begin{thm}[Multiplicativity of the Postnikov/Whitehead spectral
sequence]
\label{th:postnikov}
For cofibrant spaces $X$ and $Y$ there is a pairing of Whitehead
spectral sequences in which the $E_1$-term $E_1^{p,q} \tens E_1^{s,t}
\ra E_1^{p+s,q+t}$ is globally isomorphic to the pairing
\[
H^{q-p}(X;\pi_q E) \tens H^{t-s}(Y;\pi_t F) \ra
H^{q+t-p-s}(X\times Y;\pi_{q+t}G)
\]
up to a sign of $(-1)^{t(p-q)}$.
\end{thm}

\begin{remark}
At first glance the sign given here doesn't agree with the sign we
obtained in Theorem~\ref{th:cellular}: if we were to re-index the
Atiyah-Hirzebruch spectral sequence in the above form, the sign would
be $(-1)^{q(t-s)}$.  While this is not the same as the above sign, the
two are consistent.  Using the family of isomorphisms $(-1)^{pq}\colon
H^{q-p}(X;\pi_q E) \iso H^{q-p}(X;\pi_q E)$ and similarly for the $F$
and $G$ terms, the sign $(-1)^{t(p-q)}$ transforms into
$(-1)^{q(t-s)}$.  See Exercise~\ref{ex:signs2}.
\end{remark}

\vf 


\section{Bockstein spectral sequences}

In this section we consider two spectral sequences: one is the
classical Bockstein spectral sequence for the homotopy cofiber
sequence $H\Z \llra{\times n} H\Z \ra H\Z/n$.  The other is the
Bockstein spectral sequence for inverting the Bott element in
connective $K$-theory.

\medskip

\subsection{The Bockstein spectral sequence for $H\Z$}
Consider the following tower $(W_*,B_*)_{*\geq 0}$:

\[ \xymatrix{
&  H\Z/n & H\Z/n & H\Z/n \\
\cdots \ar[r]^n & H\Z \ar[r]^n\ar[u] & H\Z \ar[r]^n\ar[u] & H\Z. \ar[u]
}
\]
We extend this to negative degrees by taking $W_{q}=H\Z$, $B_q=*$,
and $W_{q+1} \ra W_{q}$ to be the identity map.  This is a tower of
rigid homotopy cofiber sequences, and there is an obvious pairing $(W,B)
\Smash (W,B) \ra (W,B)$ which comes from the multiplications on $H\Z$
and $H\Z/n$ (cf. \cite[Appendix C.7]{multa}).  

For any cofibrant space $X$, let $W\pair X$ denote the augmented tower
whose levels are $\F(X_+,W_{n+1}) \ra \F(X_+,W_n) \ra \F(X_+,B_n)$;
these are rigid homotopy cofiber sequences, since $H\Z$ and $H\Z/n$
are fibrant.  The homotopy spectral sequence for $W\pair X$ is called
the \mdfn{mod $n$ Bockstein spectral sequence}, and has the form
\[ E_1^{p,q}=H^{-p}(X;\Z/n) \Rightarrow H^{-p}(X;\Z).
\]
The $d_1$-differential is the usual Bockstein homomorphism.  The
multiplication on $(W,B)$ gives rise to pairings of towers $ W\pair X
\Smash W\pair Y \ra W\pair(X\times Y)$, and therefore to pairings of
spectral sequences by \cite[Thm 6.1]{multa}.  The following result is
immediate, and unlike the examples in sections 3 and 4 there are no extra
signs floating around.

\begin{thm}[Multiplicativity of the Bockstein spectral sequence]
For cofibrant spaces $X$ and $Y$ there is a pairing of Bockstein
spectral sequences whose $E_1$-term is isomorphic to the
usual pairing
$H^{-p}(X;\Z/n) \tens H^{-s}(Y;\Z/n) \ra H^{-p-s}(X\times Y;\Z/n)$ of
singular cohomology groups.  The spectral sequence converges to
the usual pairing on $H^*(\blank;\Z)$.
\end{thm}
\noindent

\subsection{The Bockstein spectral sequence for $bu$}

Let $bu$ denote a commutative ring spectrum representing connective
$K$-theory, and assume we have a map of ring spectra $bu\ra H\Z$.
Assume there is a map $\sS^2 \ra bu$ which represents a generator in
$\Ho(\sS^2,bu)\cong \Z$ (this is automatic if $bu$ is a fibrant
spectrum).  Consider the induced map
\[ \beta\colon 
\sS^2\Smash bu \llra{\beta\Smash \id} bu\Smash bu \llra{\mu} bu.
\]
It can be shown that $\sS^2\Smash bu \ra bu \ra H\Z$ is a homotopy
cofiber sequence.
If we let $(W,B)$ be the tower
\[ \xymatrixcolsep{2.5pc}
\xymatrix{
& \sS^4\Smash H\Z & \sS^2\Smash H\Z & H\Z & \sS^{-2}\Smash H\Z \\
\cdots \ar[r]^{\sS^4\Smash \beta} & \sS^4\Smash bu \ar[r]^{\sS^2\Smash
\beta}\ar[u] & \sS^2\Smash bu \ar[r]^\beta \ar[u] & bu
\ar[r]^-{\sS^{-2}\Smash \beta}\ar[u] & \sS^{-2}\Smash bu
\ar[r]^-{\sS^{-4}\Smash \beta}\ar[u]\ar[r] &\cdots }
\]
then one sees that there is a pairing $(W,B)\Smash (W,B) \ra (W,B)$
(this uses that $bu$ is commutative).  Unfortunately we are not yet in
a position to apply \cite[Thm 6.1]{multa}: $(W,B)$ is not a rigid
tower, because we don't know that $\sS^2\Smash bu \ra bu \ra H\Z$ is
null rather than just null-homotopic.  We don't get a long exact
sequence on homotopy groups until we choose null-homotopies for each
layer, and these must be accounted for.  In this particular case any
two null-homotopies are themselves homotopic, and so there should be
no problems with compatibility, but it's awkward to formulate results
along these lines.  The best way I know to proceed is actually to
discard the $H\Z$'s and consider just the tower $W_*$ consisting of
suspensions of $bu$.  We are then in a position to apply \cite[Thm
6.2]{multa}, but we must work a little harder to identify what's
happening on the $E_1$-terms.

First, one can replace $bu$ by a cofibrant commutative ring spectrum,
and this will be cofibrant as a symmetric spectrum
\cite[Thm. 4.1(3)]{SS}.  Because of this, we will just assume that our
$bu$ was cofibrant in the first place.  By \cite[B.4]{multa} there is
an augmented tower $TW_*$ consisting of cofibrations between
cofibrant spectra, a weak equivalence $TW_*\ra W_*$, and a
pairing $TW\Smash TW \ra TW$ which lifts the pairing
$W\Smash W\ra W$.  We let $C_n$ denote the cofiber of $TW_{n+1}\ra
TW_n$ and note that our pairing $TW\Smash TW \ra TW$ extends to
$(TW,C)\Smash (TW,C) \ra (TW,C)$.  

For any cofibrant space $X$, we consider the derived tower of function
spectra $\F_{der}(X_+,TW_\perp)$ from \cite[B.7]{multa}.  By
\cite[Theorems 6.1,6.10]{multa} there is an induced pairing of spectral
sequences.  Note that we have $E_1^{p,q}=\pi_p \dF(X_+,C_q) \iso \pi_p
\F(X,\Sigma^{2q}H\Z) \iso H^{2q-p}(X;\Z)$.  As usual, what we want is
to identify the pairing on $E_1$-terms with a pairing on singular
cohomology.  The following does this:

\begin{lemma}
In $\Ho(\Spectra)$ it is possible to choose a collection of
isomorphisms $C_n \ra \sS^{2n} \dSmash H\Z$ such that the following
diagrams commute:
\[ \xymatrix{
C_m \dSmash C_n \ar[rr] \ar[d] && C_{m+n}\ar[d] \\ \sS^{2m}\dSmash H\Z
\dSmash \sS^{2n}\dSmash H\Z \ar[r]^{1\dSmash t\dSmash 1} & \sS^{2m}\dSmash \sS^{2n}
\dSmash H\Z \dSmash H\Z \ar[r]^-{1\dSmash \mu} & \sS^{2(m+n)}\dSmash H\Z.  }
\]
\end{lemma}

\begin{proof}
This is very similar to the proof of Lemma~\ref{le:EMunique}, making
use of the fact that the graded ring $\oplus_n \ho(\sS^{2n},C_n)$ is
isomorphic to the ring of Laurent polynomials $\Z[t,t^{-1}]$.
Details are left to the reader.
\end{proof}

The spectral sequence  converges to
\[ \colim \bigl [ \pi_*\dF(X_+,bu) \llra{\beta\cdot} \pi_* \dF(X_+,bu)
\llra{\beta\cdot} \cdots \bigr ] 
\]
which is $\beta^{-1}[\pi_*\dF(X_+,bu)]=\beta^{-1}bu^*(X)$.  If we write
$\beta^{-1}bu^p(X)$ for the $p$th graded piece of $\beta^{-1}(bu^*X)$
then we can state the final result as follows:

\begin{thm}
For any space $X$ there is a conditionally convergent spectral
sequence
\[ E_1^{p,q}=H^{2q-p}(X;\Z) \Rightarrow \beta^{-1}bu^{-p}(X).
\]
If $X$ and $Y$ are two cofibrant spaces there is a pairing of
spectral sequences whose $E_1$-term is globally isomorphic to the
usual pairing on singular cohomology, and which converges to the
pairing $\beta^{-1}bu^*(X) \tens \beta^{-1}bu^*(Y) \ra
\beta^{-1}bu^*(X\times Y)$.  
\end{thm}

The above spectral sequence of course has the same form (up to
re-indexing) as the Atiyah-Hirzebruch spectral sequence for complex
$K$-theory.  One finds in the end that $\beta^{-1}bu^*(X)\iso K^*(X)$.  

\vf


\section{The homotopy-fixed-point spectral sequence}

Suppose $\cE$ is a fibrant spectrum with an action of a discrete group
$G$.  The homotopy fixed set $\cE^{hG}$ is defined to be
$\F(EG_+,\cE)^G$.

\begin{thm} If $\cE$ is a fibrant ring spectrum and $G$ is a discrete 
group acting via ring automorphisms, then there is a multiplicative
spectral sequence of the form $E_2^{p,q}=H^q(G;\pi_{p+q}\cE)
\Rightarrow \pi_{p}(\cE^{hG})$.  Here the pairing on $E_2$-terms is
the pairing on group cohomology with graded coefficients, defined
analagously to what was done in section~\ref{se:prelim}.
\end{thm}

\begin{remark}
The above means that the pairing on $E_2$-terms $H^q(G;\pi_{p+q}\cE)
\tens H^t(G;\pi_{s+t}\cE) \ra H^{q+t}(G;\pi_{p+q+s+t}\cE)$ is
$(-1)^{t(p+q)}$ times the `standard' pairing on group cohomology induced
by $\pi_{*}\cE \tens \pi_{*}\cE \ra \pi_{*}\cE$.
\end{remark}

\begin{proof}
Take any model for $EG$ which is an equivariant $G$-CW-complex (with $G$
acting freely on the set of cells in each dimension).  Write $EG_k$ for
$\sk_k EG$, and let $\W_n=\F(EG/EG_{n-1},\cE)$,
$\B_n=\F(EG_n/EG_{n-1},\cE)$.  The spectral sequence for this tower is
just the Atiyah-Hirzebruch spectral sequence based on $EG$ and $\cE$,
from section~\ref{se:ah}.  The sequences $\W_{n+1} \ra \W_n \ra \B_n$
are actually $G$-equivariant homotopy fiber sequences, which implies
that they give homotopy fiber sequences on $H$-fixed sets for any $H$.
We'll consider the associated tower $(\W_*^G,\B_*^G)$.

The $E_1$-term of the spectral sequence is $E_1^{p,q}=\pi_p
[\F(EG_q/EG_{q-1},\cE)^G]$.  The map
\[ \pi_p [\F(EG_q/EG_{q-1},\cE)^G] \ra \pi_p \F(EG_q/EG_{q-1},\cE) =
\ho(S^p \Smash [EG_q/EG_{q-1}],\cE) 
\] 
has its image in the $G$-fixed set of the right-hand-side, and one can
check that this gives an isomorphism $E_1^{p,q} \iso \ho(S^p\Smash
[EG_q/EG_{q-1}],\cE)^G$.  Since $EG_q/EG_{q-1}$ is a wedge of
$q$-spheres indexed by the set of $q$-cells, with $G$ action induced
by that on the indexing set, one obtains a natural isomorphism with
the group $\Hom_{G}(C_*(EG),\pi_{p+q}\cE)$.  The identification works
the same as in section~\ref{se:ah}, and the description of the
differential carries over as well.  In fact the map of towers
$(\W_*^G,\B_*^G) \ra (\W_*,\B_*)$ lets us compare the differential
with the one on the Atiyah-Hirzebruch spectral sequence, which we have
already analyzed.  So the identification of the $E_2$-term follows.

The pairing of augmented towers $\W(EG,\cE)\Smash \W(EG,\cE) \ra
\W(EG\times EG,\cE)$ from section~\ref{se:cell} obviously restricts to
fixed sets, giving $\W(EG,\cE)^G \Smash \W(EG,\cE)^G \ra \W(EG\times
EG,\cE)^G$.  This uses that $\cE\Smash \cE \ra \cE$ is
$G$-equivariant.  If one accepts that the diagonal map $EG\ra EG\times
EG$ is homotopic to a map $\Delta'$ which is both cellular and
$G$-equivariant, then we get an equivariant map of towers $\W(EG\times
EG,\cE) \ra \W(EG,\cE)$ and can restrict to fixed sets.  The existence
of $\Delta'$ follows from the $G$-equivariant
cellular approximation theorem.  Therefore we have an induced pairing
of augmented towers $\W_\perp^G\Smash \W_\perp^G \ra \W_\perp^G$.  The
identification of the product on $E_2$-terms again can be done by
comparing with the Atiyah-Hirzebruch spectral sequence for the tower
$(\W_n,\B_n)$.  Everything follows exactly as in
section~\ref{se:ah}.
\end{proof}

\begin{remark}
Instead of filtering $EG$ by skeleta, the spectral sequence can also
be constructed via a Postnikov tower on $\cE$, just as in
section~\ref{se:postnikov} (by functoriality, the Postnikov sections
of $\cE$ inherit $G$-actions).  Instead of needing a $G$-equivariant
cellular approximation for the diagonal map, one instead needs to
carry out $G$-equivariant obstruction theory.  
\end{remark}

\vf

\section{Spectral sequences from open coverings}
\label{se:open}

In this section we give a second treatment of the Atiyah-Hirzebruch
spectral sequence, together with a similar approach to the Leray-Serre
spectral sequence.  Our towers are obtained by using open coverings
and their generalization to {\it hypercovers\/}.  We will assume a
basic knowledge of hypercovers, for which the reader can consult
\cite{DI}.  Basically, one starts with an open cover $\{U_a\}$ of a
space $X$ and then chooses another open cover for every double
intersection $U_a\cap U_b$; for every resulting `triple intersection'
another covering is chosen, and so on.  All of this data is compiled
into a simplicial space $U_*$, called a hypercover.

The discussion in this section is much sloppier than the previous
ones, and should probably be improved at some point...  
\medskip

\subsection{The descent spectral sequence}
Given a hypercover $U_*$ of a space $X$ and a sheaf of abelian groups
$F$ on $X$, we let $F(U_*)$ denote the cochain complex one gets by
applying $F$ to the open sets in $U_*$.  If all the pieces
of the hypercover $U_n$ are such that $H^k_\shf(U_n,F)=0$ for $k>0$,
then Verdier's Hypercovering Theorem gives an isomorphism $H^p(F(U_*))
\ra H^p_{\shf}(X,F)$ which is functorial for $X$, $U$, and $F$.  It is
easy to explain how to get this.  First, choose a functorial,
injective resolution $I_*$ for $F$, and look at the double complex of
sections $D_{mn}=I_n(U_m)$.  This double complex has two edge
homomorphisms $I_*(X) \ra \Tot D_{**}$ and $F(U_*) \ra \Tot D_{**}$;
the two spectral sequences for the homology of a double complex
immediately show that these maps give isomorphisms on homology.
Composing them appropriately gives our natural isomorphism
$H^p(F(U_*)) \ra H^p(I_*(X))=H^p_\shf(X,F)$.

Given a fibrant spectrum $\cE$ and a hypercover $U_*$ of a $X$,
we use the simplicial space $U_*$ to set up a tower as in
section~\ref{se:simpspace}.  We let $W_n=\F(|U|/\sk_{n-1}|U|,\cE)$,
$B_n=\F(\sk_n|U|/sk_{n-1}|U|,\cE)$, and denote the resulting spectral
sequence by $E_*(U,\cE)$.  It is a theorem of \cite{DI} that the
natural map $|U_*| \ra X$ is a weak equivalence, and so the spectral
sequence converges to $\cE^*(X)$.  We'll call this the \dfn{descent
spectral sequence} based on $U$ and $\cE$.

Note that the spectral sequence is functorial in several ways.  It is
clearly functorial in $\cE$, and if $V_* \ra U_*$ is a map of
hypercovers of $X$ then there is a natural map $E_*(U,\cE) \ra
E_*(V,\cE)$.  Also, if $f\colon Y\ra X$ is a map then there is a
pullback hypercover $f^{-1}U_* \ra Y$, and a map of spectral sequences
$ E_*(U,\cE) \ra E_*(f^{-1}U,\cE)$.

So far we have said nothing about the $E_2$-term.  A space $X$ is
\dfn{locally contractible} if given any point $x$ and any open
neighborhood $x\in V$, there is an open neighborhood $x\in W \subseteq
V$ such that $W$ is contractible.  Given such a space $X$, one can
build a hypercover $U_*$ of $X$ in which every level is a disjoint
union of contractible opens.  We'll call $U_*$ a
\dfn{contractible hypercover}.  

If $U_* \ra X$ is a hypercover there is a natural isomorphism
$E_2^{p,q}(U,\cE)\cong H^q({{\cE^{-p-q}}}(U_*))$.  If we assume $X$ is
locally contractible and $U_*$ is a contractible hypercover, we can
simplify this further.  Let ${\tcE^q}$ denote the sheafification of
the presheaf $V \mapsto \cE^q(V)$, and note that this is a locally
constant sheaf on $X$.  It follows that $\cE^q(V) \ra {\tcE^q}(V)$ is
an isomorphism for every contractible $V$.  In particular,
$H^q(\cE^{p-q}(U_*)) \ra H^q(\tcE^{p-q}(U_*))$ is an isomorphism.  But
since the sheaves $\tcE^*$ are locally constant, it follows that if
$V$ is a contractible open set in $X$ then $H^p_\shf(V,\tcE^q)=0$ for
$p>0$ (see \cite[II.11.12]{Br}).  The fact that $U_*$ is a
contractible hypercover now shows that 
\[ E_2^{p,q}(U,\cE) \iso H^q(\cE^{p-q}(U_*)) \iso H^q(\tcE^{p-q}(U_*)) \iso
H^q_{\shf}(X,{{\tcE}^{-p-q}})
\]
(and all these isomorphisms are natural). 

\begin{remark}
\label{re:HC}
Let us for a moment forget about contractible hypercovers, and look at
all of them.  For each map of hypercovers $V_*\ra U_*$ there is an
induced map of spectral sequences $E_*(U,\cE) \ra E_*(V,\cE)$.  Two
homotopic maps $V_* \dbra U_*$ induce the same map of spectral
sequences from $E_2$ on.  So if we forget about $E_1$-terms, we have a
diagram of spectral sequences indexed by the homotopy category of
hypercovers $\HC_X$.  We can take the colimit of these spectral
sequences, and the Verdier Hypercovering Theorem identifies the
$E_2$-term with sheaf cohomology.  So the hypothesis that the base $X$
be locally contractible is not really necessary for the development of
our spectral sequence.  It does make things easier to think about,
though, so we will continue to specialize to that case.
\end{remark}

\medskip

Now suppose that $\cE$ is a fibrant ring spectrum, $Y$ is another
locally contractible space, and $V_* \ra Y$ is a contractible
hypercovering.  The product simplicial space $U_*\times V_*$ is a
contractible hypercover of $X\times Y$.  By the discussion in
section~\ref{se:simpspace} we get a pairing of spectral sequences
$E_*(U,\cE) \tens E_*(V,\cE) \ra E_*(|U|\times |V|,\cE)$ and a map of
spectral sequences $E_*(U\times V,\cE) \ra E_*(|U|\times |V|,\cE)$.  
It follows from Proposition~\ref{pr:simp} that the latter map is an
isomorphism (from the hypercover $U$ one defines a simplicial set
$\pi_0U$ by applying $\pi_0(\blank)$ in every dimension, and in our
case $U\ra \pi_0(U)$ is a levelwise weak equivalence; the same holds
for $V$ and $U\times V$, and so we are really just dealing with
simplicial sets).  Therefore we get the following:

\begin{thm}[Multiplicativity of the descent spectral sequence]
\label{th:descent}
Let $X$ and $Y$ be locally contractible spaces, with contractible
hypercovers $U_*\ra X$ and $V_*\ra Y$.  Then there is a
pairing of descent spectral sequences $E_*(U,\cE)\tens E_*(V,\cE) \ra
E_*(U\times V,\cE)$ in which the $E_2$-term is globally isomorphic
to
\[ H^q_{\shf}(X,\cE^{-p-q}) \tens H^t_{\shf}(Y,\cE^{-s-t}) \ra
H^{q+t}_{\shf}(X\times Y, \cE^{-p-q-s-t})
\]
up to a sign difference of $(-1)^{t(p+q)}$.
\end{thm}  

\begin{proof}
To identify the product on $E_2$-terms we use the fact that if $F$ and
$G$ are two sheaves on $X$ and $Y$, and $U$ and $V$ contractible
hypercovers of $X$ and $Y$, then the pairing of cosimplicial abelian
groups $F(U_*)\tens G(V_*) \ra (F\tens G)(U_* \times V_*)$ induces
the cup product on sheaf cohomology via an Eilenberg-Zilber map
and the isomorphisms $H^p(F(U_*))\iso H^p(X,F)$, etc.  Someone should
write down a careful proof of this someday, but I have proven enough
`trivial' things for one paper.  

The sign comes from an implicit use of the suspension isomorphism.
The $E_2$-term of $E_*(U,\cE)$ is most naturally identified with
$H^q_{\shf}(X,\cG^{p,q})$ where $\cG^{p,q}$ denotes the sheafification
of $V\mapsto \cE^{-p}(S^q\Smash V_+)$.  Clearly $\cG^{p,q}$ is
isomorphic to the sheaf $\tcE^{-p-q}$ via the left-suspension
isomorphism, but this introduces signs into the pairings much like in
Theorem~\ref{th:postnikov}.  We leave the details to the reader.
\end{proof}

The assumption in the above theorem that $X$ and $Y$ be locally
contractible is probably not really necessary, as in
Remark~\ref{re:HC}; but I have not worked out the details.  Also, when
$X$ and $Y$ are both locally contractible and paracompact, sheaf
cohomology with locally constant coefficients is isomorphic to
singular cohomology; the pairing on sheaf cohomology corresponds to
the cup product under this isomorphism \cite[Theorem III.1.1]{Br}.  So
the descent spectral sequence takes the same form as the
Atiyah-Hirzebruch spectral sequence in this case.

\begin{exercise}
If $U_* \ra X$ is a contractible hypercover, let $\pi_0(U)$ be the
simplicial set defined above---obtained by replacing each $U_n$ by its
set of path components.  There is a natural map $U\ra \pi_0(U)$, and
this is a levelwise weak equivalence.  Convince yourself that the
descent spectral sequence $E_*(U,\cE)$ is canonically isomorphic to
the Atiyah-Hirzebruch spectral sequence for the space $|\pi_0(U)|$
based on its skeletal filtration.  Conclude that everything in this
section is just a restatement of the results in section \ref{se:ah}.
\end{exercise}

\subsection{The Leray spectral sequence}

Now assume that $\pi\colon E\ra B$ is a fibration with fiber
$F$, and that $B$ is locally contractible.  Choose a hypercover $U_*$
of $B$ and consider the pullback hypercover $\pi^{-1}U_*$ of $E$.
Once again, we know by \cite{DI} that $|\pi^{-1}U_*| \ra E$ is a weak
equivalence.  Taking $\cE=H\Z$, the skeletal filtration of
$|\pi^{-1}U_*|$ gives a spectral sequence for computing the
homotopy groups of the function space $\F(E,H\Z)$---i.e., it computes
the singular cohomology groups of $E$.  This is the \dfn{Leray (or
Leray-Serre) spectral sequence}.

\begin{exercise}
Let $H^n\pi^{-1}$ denote the sheaf on $B$ obtained by sheafifying
$V\mapsto H^n_{\sng}(\pi^{-1}V)$.  Check that $H^n\pi^{-1}$ is a
locally constant sheaf on $B$ whose stalks are isomorphic to
$H^n_{\sng}(F)$.  So if $B$ is simply-connected, it is a constant
sheaf.  We will abbreviate $H^n\pi^{-1}$ as $\cH^n(F)$.  
\end{exercise}

It follows that if $U$ is a contractible hypercover then the
$E_2$-term of the spectral sequence is isomorphic to
$H^q_\shf(B,\cH^{-p-q}(F))$.  When $B$ is simply connected this is
$H^q_{shf}(B,H^{-p-q}F)$, and when $B$ is paracompact this can be
identified with singular cohomology.

Once again, if $E' \ra B'$ is a second fiber bundle with fiber $F'$,
there is a pairing of spectral sequences which on $E_2$-terms has the
form
\[ H^q_\shf(B,\cH^{-p-q}F) \tens H^t_\shf(B',\cH^{-s-t}F') \ra
H^{q+t}_\shf(B\times B',\cH^{-p-q-s-t}(F\times F')).
\]
One has the usual sign difference from the canonical pairing, for the
same reasons as in Theorem~\ref{th:descent}.  If $E'\ra B'$ is the
same as $E\ra B$, then we can compose with the diagonal map to get a
multiplicative structure on the Leray-Serre spectral sequence for
$E\ra B$.

\begin{thm}[Multiplicativity of Leray-Serre]
\label{th:leray}
Let $E\ra B$ and $E'\ra B'$ be fibrations where $B$ and $B'$ are
locally contractible.  Then there
is a pairing of Leray-Serre spectral sequences for which the
$E_2$-term is globally isomorphic to the pairing
\[ H^q_{\shf}(B,\cH^{-p-q}F) \tens H^t_{\shf}(B',\cH^{-s-t}F') \ra
H^{q+t}_{\shf}(B\times B',\cH^{-p-q-s-t}(F\times F')).
\]
except for  a sign difference of $(-1)^{t(p+q)}$.   
\end{thm}

\begin{exercise}
Fill in the many missing details from this section.
\end{exercise}

\begin{exercise}
Of course we didn't need to take $\cE=H\Z$ in the above discussion.
We could have used any ring spectrum, in which case we would obtain
the combination Atiyah-Hirzebruch-Leray-Serre spectral sequence.
Think through the details of this.
\end{exercise}

\vf


\bibliographystyle{amsalpha}

\end{document}